\documentclass[aos,MSNbibl,nameyear,dvips]{arximspdf}
\usepackage{graphicx}

\doi{10.1214/12-AOS1074} 
\volume{41}
\issue{1}
\pubyear{2013}
\firstpage{63}
\lastpage{90}

\makeatletter

\def\cal{\mathcal}
\def\bm{\bolds}
\def\open{}
\def\E{{\open\mathbb{E}}}
\def\proj{\operatorname{pr}^{\bot}}

\newtheorem{theorem}{Theorem}
\newproclaim{remark}{Remark}
\newtheorem{lemma}{Lemma}
\newproclaim{definition}{Definition}
\newtheorem{corol}{Corollary}
\newproclaim{assumption}{Assumption}

\makeatother

\begin{document}
\begin{frontmatter}

\title{Universally optimal crossover designs under subject dropout}
\runtitle{Crossover under dropout}

\begin{aug}
\author{\fnms{Wei} \snm{Zheng}\corref{}\ead[label=e1]{weizheng@iupui.edu}}
\runauthor{W. Zheng}
\affiliation{Indiana University-Purdue University Indianapolis}
\address{Department of Mathematical Sciences\\
Indiana University-Purdue University Indianapolis\\
Indianapolis, Indiana 46202-3216\\
USA\\
\printead{e1}}
\end{aug}

\received{\smonth{4} \syear{2012}}
\revised{\smonth{9} \syear{2012}}

%
\begin{abstract}
Subject dropout is very common in practical applications of cross\-over
designs. However, there is very limited design literature taking this
into account. Optimality results have not yet been well established due
to the complexity of the problem. This paper establishes feasible, as
well as necessary and sufficient conditions for a crossover design to
be universally optimal in approximate design theory in the presence of
subject dropout. These conditions are essentially linear equations with
respect to proportions of all possible treatment sequences being
applied to subjects and hence they can be easily solved. A general
algorithm is proposed to derive exact designs which are shown to be
efficient and robust.
\end{abstract}

%


\end{frontmatter}

\section{Introduction}\label{sec1}

Crossover designs have been widely used in industry due to their cost
effectiveness and statistical efficiency. They are applicable for
experiments aiming to compare effects of different treatments by
applying them to a number of subjects across several periods. The
response observation is typically modeled by additive effects of
subjects, periods, treatments and the carryover effects of the
treatment from the previous period. There has been tremendous amount of
literature regarding the identification of optimal designs. See \citet{HedAfs78},
Ch{\^e}ng and Wu (\citeyear{CheWu80}),
\citet{Kun84}, \citet{Stu91}, Kushner (\citeyear{Kus97N1}, \citeyear{Kus97N2}, \citeyear{Kus98}),
\citet{KunMar00}, \citet{KunStu02}, Hedayat and Yang (\citeyear{HedYan03}, \citeyear{HedYan04}) and
\citet{HedZhe10}, for instance. For comprehensive reviews, see \citet{Mat88},
\citet{RatEvaAll92}, \citet{Stu96}, Jones and
Kenward (\citeyear{JonKen89}), \citet{Sen03} and \citet{BosDey09}.

An important issue regarding crossover designs is that subject may drop
out of the study. As a result, the experiment will not be carried out
as planned. \citet{Mat88} commented this is one of the main concerns
of crossover designs. \citet{LowLewPre99} observed that ``A
dropout rate of between $5\%$ and $10\%$ is not uncommon and, in some
areas, can be as high as $25\%$.'' Meanwhile, a design, which is optimal
or highly efficient in the absence of dropout, would become inefficient
or even disconnected in the presence of subject dropout. Examples could
be found in \citet{God04}, Majumdar, Dean and Lewis (\citeyear{MajDeaLew08}) as well
as Section~\ref{sec6} of this paper.

To conclude, it is very important to find optimal or efficient designs
in the presence of subject dropout, yet there is very limited
literature on this. \citet{BosBag08} derived designs which are
universally optimal for both direct and carryover effects for both the
situation of no dropout and the situation that all subjects drop out
after period $q$ with $q$ being judiciously chosen. Similar results are
presented by Majumdar, Dean and Lewis (\citeyear{MajDeaLew08}). The latter restricted the
comparison of designs within the subclass of uniformly balanced
repeated measurement designs (UBRMDs), whose optimality property has
been well recognized in literature for the situation of no dropout. For
the second situation with any given $q$, they proposed type ${\cal
W}_q$ UBRMDs, which reduce the maximum loss of the information for
parameters in terms of $A$-criterion as compared to general UBRMDs.
Following the latter paper, Zhao and Majumdar (\citeyear{ZM12}) further explored
the special case when~$q$ is one less the number of periods and the
numbers of treatments and periods are the same.

The previous three papers share two drawbacks: (i) The proposed
designs exist only under very rare combinations of the numbers of
subjects, periods and treatments. See Section~\ref{sec41} for
relevant discussions of the former paper. As for the other two papers,
it is well known that the existence of UBRMDs is rare. (ii) The
information regarding the mechanism of how subjects drop out was not
taken into account.

To address the latter drawback, it is plausible to measure the
performance of designs by taking the expectation of a regular
optimality criterion with respect to the mechanism of subject dropout.
\citet{LowLewPre99} worked in this direction by using
intensive computer programming. They concluded that when the Latin
squares consisting of the design is more diverse, the resulting design
performs better in terms of both efficiency and robustness. This
argument is further supported by the comparison in Section~\ref
{sec6}. However, the case studies they provided fail to provide
general guidance in identifying efficient designs. To serve this
purpose, theoretical results are called for.

In this paper, we develop feasible equivalent conditions for a design
to be universally optimal for direct treatment effects in approximate
design theory under the same setup as that of \citet{LowLewPre99}.
The equivalence holds for any probability distribution of
subject dropout. The results can be easily modified to find optimal or
highly efficient exact designs for any combination of the numbers of
subjects, periods and treatments. As a result, the two drawbacks are
both addressed here.\vadjust{\goodbreak}

The rest of the paper is organized as follows. Section~\ref{sec2}
formulates the problem, introduces notation and gives some preliminary
results. Section~\ref{sec3} introduces necessary concepts in
approximate design theory, proves the existence of universally optimal
designs and also gives necessary, sufficient and equivalent conditions
for universal optimality. Section~\ref{sec5} gives explicit and
feasible forms of optimality conditions in terms of linear equations,
which are built upon the preceding section. Section~\ref{sec6}
further provides a general algorithm for deriving an optimal or
efficient exact design for any combination of the numbers of subjects,
periods and treatments as well as any probability distribution of
subject dropout. Besides, comparisons are made to designs in
literature. Section~\ref{sec7} summarizes the results. Finally, some
proofs are deferred to Section~\ref{sec8}.

\section{Framework}\label{sec2}

This section introduces the framework of the problem. Section \ref
{sec21} introduces the statistical model for the design problem and
provides notation and assumptions necessary to the rest of the paper.
Section~\ref{sec22} defines an ideal target function in finding a
design, proposes a corresponding surrogate target function, and
discusses the relationship between these two target functions. Section
\ref{sec23} provides some preliminary results as a preparation for
the rest of the paper.

\subsection{Modeling and notation}\label{sec21}

In a crossover design with $p$ periods, $t$ treatments and $n$
subjects, the response is typically modeled as
%
\begin{equation}
\label{model} Y_{dku}=\mu+\pi_k+\varsigma_u+
\tau_{d(k,u)}+\gamma_{d(k-1,u)}+\varepsilon_{ku},
\end{equation}
where $\{\varepsilon_{ku}, 1\leq k\leq p, 1\leq u\leq n\}$ are
independent with mean zero and variance~$\sigma^2$. Here, $Y_{dku}$
denotes the response from subject $u$ in period $k$ to which treatment
$d(k,u)\in\{ 1,2,\ldots,t \}$ was assigned by design $d$. Furthermore,
$\mu$ is the general mean,
$\pi_k$ is the $k$th period effect, $\varsigma_u$ is the $u$th
subject effect, $\tau_{d(k,u)}$ is the (direct) treatment effect of
treatment $d(k,u)$ and $\gamma_{d(k-1,u)}$ is the
carryover effect of treatment $d(k-1,u)$ that subject $u$ received in
the previous period (by convention $\gamma_{d(0,u)}=0$).

Let $G$ be a temporary object whose meaning differs from context to
context. Then we define $G'$ to represent the transpose of the matrix~$G$,
$G^-$ to represent a generalized inverse of the matrix $G$,
$\operatorname{ tr}(G)$ to represent the trace of the matrix~$G$ and $\proj$ to be a
projection operator such that $\proj G=I-G(G'G)^-G'$. For two square
matrices of equal size, $G_1$ and $G_2$, $G_1\leq G_2$ means that
$G_2-G_1$ is nonnegative definite. For a set $G$, the number of
elements in the set is represented by~$|G|$.

Besides, $I_k$ is the $k\times k$ identity matrix, $1_k$ is the vector
of length $k$ with all its entries as $1$, $J_k=1_k1_k'$ is the square
matrix with all its entries as $1$. We further define $B_k=I_k-J_k/k$,
$B_{ij}^k$ to be the $i\times j$ matrix with its upper left corner
filled with the submatrix $B_k$ while the remaining entries filled with
$0$, and $B^k_i=B^k_{ii}$. The notation of $I_{ij}^k$ and $I_i^k$ are
defined in the same fashions as $B_{ij}^k$ and $B_i^k$. Finally,
$\otimes$ represents the Kronecker product of two matrices. To make
the problem resolvable, it is necessary to make two mild assumptions as follows.

\begin{assumption}\label{ass1}
Once a subject drops out of the study, the probability that the subject
reenters the study is zero.
\end{assumption}

By Assumption~\ref{ass1}, we are able to define $l_i$, $1\leq i \leq
n$, to be the total number of periods that subject $i$ stayed in the
experiment. Further it is realistic in a large number of applications
to assume the following:

\begin{assumption}\label{ass2}
The dropping out mechanism is independent of the choice of design $d$
as well as the outcome of the experiments. Moreover $\{l_i,1\leq i\leq
n\}$ are i.i.d.
\end{assumption}

By Assumption~\ref{ass2}, we could define $a_k$ to be the probability
that $l_i=k$, $1\leq k\leq p$, and hence we are in place to define the
following technical terms:

\begin{itemize}
\item$\vec{a}=(a_1,a_2,\ldots,a_p)$.
\item$a_{jk}=\sum^k_{i=j}a_i$, $1\leq j\leq k\leq p$. (Convention:
$a_{p+1,p}=0$.)
\item$m=\min\{k\dvtx a_k>0\}$.
%
\item$\alpha_k=n^{-1}
((n+1)a_k+a_{1,k-1}^{n+1}-a_{1k}^{n+1} )$, $1\leq k\leq p$.
\item$\beta_k=a_k+a_{k+1,p}a_{1k}^n-a_{kp}a_{1,k-1}^n$, $1\leq k\leq p$.
\item$A=\sum^p_{k=1}\alpha_kB^k_p$.
\item$B=\sum^p_{k=1}\beta_kB^k_p$.
\end{itemize}

\begin{definition}
An experiment is said to be complete if there is no dropout.
\end{definition}

By definition the complete experiment is a special case in our
framework and has been extensively studied in literature. Here, we aim
to investigate desirable designs for any given dropout mechanism $\vec{a}$.

Notice that $A$ and $B$ are both nonnegative definite matrices. Since
$\beta_k\geq a_k+a_{k+1,p}a_{1k}^n-a_{kp}a_{1k}^n=a_k(1-a_{1k}^n)\geq
0$, we have $B\geq0$. By the mean value theorem one could show that
$\alpha_k\geq0$ and hence $A\geq0$. Note that $a_k=0$ implies
$\alpha_k=\beta_k=0$. Hence we have $A=\sum^p_{k=m}\alpha_kB^k_p$
and $B=\sum^p_{k=m}\beta_kB^k_p$. The same representation will be
adopted in the sequel whenever the summation over the period $k$ is
involved. Finally, we should be aware of the differences and
relationships among the matrices $B_k$, $B^k_p$ and $B$.

\subsection{Optimality criteria}\label{sec22}
Writing the $np\times1$ response vector as ${
Y}_d=(Y_{d11}, Y_{d21},\ldots,Y_{dp1},Y_{d12},\ldots,Y_{dpn})'$, model (\ref
{model}) can be written as
%
\begin{equation}
\label{model2} Y_d ={ 1}_{np}\mu+Z\bm{\pi}+U\bm{
\varsigma}+{ T}_d\bm{\tau}+{F}_d\bm{\gamma}+
\varepsilon,\vadjust{\goodbreak}
\end{equation}
where $\bm{\pi}=(\pi_1,\ldots,\pi_p)'$, $\bm{\varsigma}=(\varsigma_1,\ldots,\varsigma_n)'$, $\bm{\tau}=(\tau_1,\ldots,\tau_t)'$, $\bm
{\gamma}=(\rho_1,\ldots,\rho_t)'$,
$Z={ 1}_n\otimes{ I}_p$, $U=I_n\otimes1_p$ and ${ T}_d$ and ${ F}_d$
denote the treatment/subject and carryover/subject incidence matrices.
Here $\E\varepsilon=0$ and $\operatorname{Var}(\varepsilon)=\sigma^2I_{np}$. For
design $d$ under a realization of experiment $l=(l_1,\ldots,l_n)'$, the
information matrix for the direct treatment effects $\tau$ under model
(\ref{model2}) with $\sigma^2=1$ is
\begin{eqnarray*}
C_d(\tau,l)&=&(MT_d)'
\proj(MZ|MU|MF_d) (MT_d)
\\
&=&C_{d11}(l)-C_{d12}(l)\bigl[C_{d22}(l)
\bigr]^-C_{d21}(l),
\end{eqnarray*}
where
\begin{eqnarray*}
C_{d11}(l)&=&T_d'OT_d, \qquad C_{d12}(l)=T_d'OF_d,\\
C_{d21}(l)&=&C_{d12}',\qquad C_{d22}(l)=F_d'OF_d,\\
M&=&\operatorname{diag}\bigl(I_{l_i,p}^{l_i},i=1,2,\ldots,n\bigr),\\
O&=& M'\proj(MZ|MU)M.
\end{eqnarray*}

Under a complete experiment, \citet{Kie75} defined a design to be
universally optimal if it maximizes $\Phi(C_d(\tau,p1_n))$ for any
$\Phi$ satisfying:
\begin{longlist}[(C.1)]
\item[(C.1)] $\Phi$ is concave;
\item[(C.2)] $\Phi(S'CS)=\Phi(C)$ for any permutation matrix $S$;
\item[(C.3)] $\Phi(bC)$ is nondecreasing in the scalar $b>0$.
\end{longlist}

Optimality criteria defined by such a $\Phi$ includes, but is not
limited to, $A$, $D$, $E$ and $T$. See \citet{Kie75} and \citet{Yeh86} for
instance. In the subject dropout setup there does not exist a design
which maximizes $\Phi(C_d(\tau,l))$ for all realizations of $l$. One
reasonable target is to find a design which maximizes $\phi_0(d|\Phi
,\vec{a}):=\E_{\vec{a}}\Phi(C_d(\tau,l))$ for any $\Phi$
satisfying the above three conditions. Here the expectation is taken
over the probability space of $l$ with parameter~$\vec{a}$. For
notational simplicity, we would omit the subscript $\vec{a}$ for $\E$
and the parameters $\Phi$ and $\vec{a}$ for $\phi_0$ whenever it is
clear from the context. So we have $\phi_0(d):=\phi_0(d|\Phi,\vec
{a})=\E\Phi(C_d(\tau,l))$.

There are two major difficulties in maximizing $\phi_0(d)$ which make
the problem intractable, if not impossible: (i) $\Phi$ is a
nonlinear function and hence the expectation would interact with the
form of $\Phi$. (ii) Even when the dropout situation $l$ is fixed,
there is still a lack of tools to deal with the information matrix
$C_d(\tau,l)$ if subjects drop out at different periods under $l$. In
order to tackle these difficulties, we propose to replace the original
target function of $\phi_0(d)$ with the surrogate target function of
$\phi_1(d)=\Phi(C_d)$ where
%
\begin{eqnarray}\label{eqn418}
C_d&=&C_{d11}-C_{d12}C_{d22}^-C_{d21},
\nonumber
\\[-8pt]
\\[-8pt]
\nonumber
C_{dij}&=&\E C_{dij}(l),\qquad 1\leq i,j\leq2.
\end{eqnarray}
It will be shown in Section~\ref{sec6} that this replacement is very
successful in identifying highly efficient, if not optimal, designs for
the criterion $\phi_0(d)$.\vadjust{\goodbreak}

For $i=0$ or $1$, let $d^*_i$ be an optimal design under $\phi_i$.
Then define $e_i(d)=\phi_i(d)/\phi_i(d^*_i),i=0,1$, to be the
efficiency of $d$ under $\phi_i$-criterion. Also we call $g(d)=\phi_0(d)/\phi_1(d)$ to
be the \textit{gap} function between the two target
functions for design $d$. Even though we are working on $\phi_1$
instead of $\phi_0$, the $\phi_0$-efficiency $e_0(d)$ could be
bounded by $e_1(d)g(d)$ as shown by Lemma~\ref{lemma203}.

\begin{lemma}[{[\citet{Puk93}, pages 74--77]}]\label{lemma2122}
The Schur complement of a matrix $G\geq0$ is a concave nondecreasing
function of $G$.
\end{lemma}

\begin{lemma}\label{lemma203}
For any $\Phi$, $\vec{a}$ and design $d$, we have $\phi_0(d)\leq
\phi_1(d)$. Further we have
\[
e_0(d)\geq e_1(d)g(d).
\]
In particular, for any $\phi_1$-optimal design $d$, we have
$e_0(d)\geq g(d)$.
\end{lemma}
\begin{pf} By Lemma~\ref{lemma2122} we have $\E C_d(\tau,l)\leq
C_d$. Then we have
%
\begin{eqnarray}
\phi_0(d)&=&\E\Phi\bigl(C_d(\tau,l)\bigr)
\nonumber
\\
&\leq& \Phi\bigl(\E C_d(\tau,l)\bigr)\label{eqn123}
\\
&\leq& \Phi(C_d)\label{eqn1232}
\\
&=&\phi_1(d).\label{eqn202}
\end{eqnarray}
By (\ref{eqn202}) we have
$e_0(d)=\phi_0(d)/\phi_0(d^*_0)
\geq\phi_0(d)/\phi_1(d^*_0)
\geq\phi_0(d)/\phi_1(d^*_1)
=\break  e_1(d)g(d).$
\end{pf}

By (\ref{eqn202}) we have $g(d)\leq1$, and hence $g(d_1^*)\leq
e_0(d^*_1)$. That means if we could find a $\phi_1$-optimal design,
then the value of the gap function $g$ evaluated at this design serves
as a lower bound of its $\phi_0$-efficiency. Inequalities (\ref
{eqn123}) and (\ref{eqn1232}) are essentially Jensen-type
inequalities. The equalities therein both hold if the realization of
subject dropout, $l$, is not random. When the variation in $l$ is not
very large, it would be plausible to work on the surrogate target of
maximizing $\phi_1(C_d)$ instead of $\phi_0(C_d)$ since the value of
the gap function $g$ would be close to unity. Note that a popular
choice of $\Phi$ is the trace of a matrix ($T$-criterion), for which the
equality in (\ref{eqn123}) always holds.

When the experiment is complete, the necessary and sufficient
conditions for $\phi_1$-universal optimality derived in Section \ref
{sec5} reduce to that of \citet{Kus97N2}. Note that the matrix $C_d$
in (\ref{eqn418}) is no longer an information matrix for any design,
and as a result the ideas of proving the existence of universally
optimal designs, given by Theorem 3.4 of \citet{Kus97N2}, are not
applicable here. However, we found that similar results could be
derived by direct manipulation on the matrix~$C_d$. See Sections \ref
{sec32} and~\ref{sec33} for details. Moreover, since $A\neq B$ in
general, the arguments in deriving the linear equation as in proof of
Theorem 5.3 of \citet{Kus97N2} are not applicable here either. For the
approach of tackling this difficulty, see Section~\ref{sec51} for details.

\subsection{Preliminary results}\label{sec23}

\begin{lemma}\label{lemma125}
Under Assumptions~\ref{ass1} and~\ref{ass2} we have
$C_{d11}=T_d'VT_d$, $C_{d12}=T_d'VF_d$ and $C_{d22}=F_d'VF_d$ with
%
\begin{eqnarray}
\label{eqn2162} V&=&\sum^p_{k=m}\bigl(
\alpha_kI_n-n^{-1}\beta_kJ_n
\bigr)\otimes B^k_p
\nonumber
\\[-8pt]
\\[-8pt]
\nonumber
&=&I_n\otimes A-n^{-1}J_n\otimes B.
\end{eqnarray}
\end{lemma}

Since $B^1_p=0$ the $m$ in (\ref{eqn2162}) could be replaced by
$\operatorname{max}(m,2)$. A heuristic explanation for this observation is that when
$l_i=1$ there is no information gained from this subject, because we
rely on within subject comparison for treatments in crossover designs.
When the experiment is complete we have $\alpha_k=\beta_k=0$ for all
$1\leq k\leq p-1$ and $\alpha_p=\beta_p=1$. In this case, we have the
reduction of $A=B=B_p$ and $V=B_n\otimes B_p$, for which the optimality
problem has been extensively studied in literature.

\begin{corol}
Any design which is $\phi_1$-optimal with $\Phi$ satisfying
conditions~\textup{(C.1)--(C.3)} under model (\ref{model}) is still optimal
under the same criterion when the within subject covariance is of the
form
%
\begin{equation}
\label{eqn21210} \Sigma=I_p+\eta1_p'+1_p
\eta'.
\end{equation}
One special case is the compound symmetric covariance matrix, that is,
$\Sigma=I_p+bJ_p$. Here $\eta$ is an arbitrary vector, and $b$ is an
arbitrary real number.
\end{corol}

\begin{pf} Let $\Sigma_k$ be the $k\times k$ upper left submatrix of
$\Sigma$ for $1\leq k\leq p$. By direct calculation, we have
%
\begin{equation}
\label{eqn125} \Sigma_k^{-1}-\Sigma_k^{-1}J_k
\Sigma_k^{-1}/1_k'
\Sigma_k^{-1}1_k= B_k.
\end{equation}
By following the same calculation as the proof of Lemma \ref
{lemma125}, the corollary is established in view of equation (\ref{eqn125}).
\end{pf}

%
\begin{remark} The covariance matrix as in (\ref{eqn21210}) is
called a ``type-H'' matrix; see \citet{HuyFel70}.
\end{remark}

\section{\texorpdfstring{$\phi_1$-universal optimality}{phi1-universal optimality}}\label{sec3}

This section explores the $\phi_1$-universal optimality in approximate
design theory, where $\phi_1$-universal optimality is defined as follows.
%
\begin{definition}
Given $p$, $t$, $n$ and a dropout mechanism $\vec{a}$, a design $d$ is
said to be $\phi_1$-universally optimal if $d$ maximizes $\phi_1(d)$
over all designs for any $\Phi$ satisfying conditions (C.1), (C.2) and (C.3).
\end{definition}

Section~\ref{sec31} introduces the ideas in approximate design theory
as well as the concept of symmetric designs. Section~\ref{sec32}
shows that a design would be $\phi_1$-universally optimal as long as
its information matrix is of the form $C_d=ny^*B_t/(t-1)$ with $y^*$
introduced by equation (\ref{eqn2037}). Section~\ref{sec33} shows
that there always exists a symmetric design which satisfies this
sufficient condition for $\phi_1$-universal optimality, and further by
argument of \citet{Kie75} that this condition is also necessary for any
design to be $\phi_1$-universally optimal. However, this condition is
not immediately applicable for application. Section~\ref{sec5} gives
an equivalent condition which is more readily applicable. Some relevant
technical preparations are given in Sections~\ref{sec34} and~\ref{sec35}.

\subsection{Approximate design theory and symmetric designs}\label{sec31}

A design $d$ with $p$ periods, $t$ treatments and $n$ subjects could be
considered as the result of selecting~$n$ sequences with replacement
from the collection of all possible $t^p$ sequences, and this
collection is denoted by ${\cal S}$. Let $n_s$ be the number of
replications of sequence $s$ in the design, and define $P_d=(p_s,s\in
{\cal S})$ with $p_s=n_s/n$. When we ignore the ordering of the $n$
sequences in the design, we have the one to one correspondence of
$d\leftrightarrow(n,P_d)$ with the restrictions of (i) $\sum_{s\in
{\cal S}}p_s=1$, (ii) $p_s\geq0$ and (iii) $np_s$ being an integer
for all $s$. In approximate design theory, we only keep the first two
restrictions and allow $np_s$ not to be an integer.

Let $\sigma$ be a permutation of symbols $\{1,2,\ldots,t\}$. For a
sequence $s=(t_1,\ldots,\break t_p)$, we define $\sigma s=(\sigma
(t_1),\ldots,\sigma(t_p))$. Then the design $\sigma d$ is defined by
$P_{\sigma d}=(p_{\sigma^{-1}s},s\in{\cal S})$. The permutation
matrix $S_{\sigma}$ is the unique matrix satisfying $T_{\sigma
s}=T_sS_{\sigma}$ for all $s\in{\cal S}$. In the sequel we replace
the subject index $u$ by sequence index $s$ whenever it is necessary.

A design $d$ is said to be symmetric if $P_d=P_{\sigma d}$. Also we
define symmetric blocks as $\langle s\rangle=\{\sigma s, \sigma\in
{\cal P}\}$ where ${\cal P}$ is the collection of all possible $t!$
permutations, that is, $|{\cal P}|=t!$. We further define $p_{\langle
s\rangle}=\sum_{\tilde{s}\in\langle s\rangle}p_{\tilde{s}}$. For
a symmetric design, we have $p_{\tilde{s}}=p_{\langle s\rangle
}/|\langle s\rangle|$ for any $\tilde{s}\in\langle s\rangle$. Given
$p,t,n$, a symmetric design $d$ is uniquely determined by $(p_{\langle
s\rangle},s\langle\in\rangle{\cal S})$, where $s\langle\in\rangle
{\cal S}$ means that $s$ runs through all distinct symmetric blocks
contained in ${\cal S}$.

\subsection{\texorpdfstring{A sufficient condition for $\phi_1$-universal optimality}
{A sufficient condition for phi1-universal optimality}}\label{sec32}

Denote by $T_u$ (resp., $F_u$) the $p\times t$ submatrix of $T$ (resp., $F$)
corresponding to the $u$th subject. Define $\overline{T}=n^{-1}\sum^n_{u=1}T_u$, $\hat{T}_u=T_uB_t$ and $\hat{\overline{T}}=\overline
{T}B_t$. The notation $\overline{F}$, $\hat{F}_u$ and $\hat
{\overline{F}}$ are defined in the same way corresponding to carryover
effects. Let $\overline{C}_{dij}=\sum_{\sigma\in{\cal P}}S_{\sigma
}'C_{dij}S_{\sigma}/|{\cal P}|$, $1\leq i,j\leq2$ and $\overline
{C}_d=\sum_{\sigma\in{\cal P}}S_{\sigma}'C_dS_{\sigma}/|{\cal
P}|$. Note that $\overline{C}_{dij}, 1\leq i,j\leq2$, are completely
symmetric, also $\overline{C}_{d11}$ and $\overline
{C}_{d12}=(\overline{C}_{d21})'$ have row and column sums as zero. Let
$\mathbb{I}$ be the indicator function. By Proposition 1 of \citet{KunMar00}, we have
%
\begin{eqnarray}
\overline{C}_d&\leq&\tilde{C}_d\label{eqn2034}
\\
&=& \biggl(c_{d11}-\frac{c_{d12}^2}{c_{d22}}\mathbb {I}_{[c_{d22}>0]} \biggr)
\frac{B_t}{t-1},\label{eqn2035}
\end{eqnarray}
where
\begin{eqnarray*}
\tilde{C}_d&=&\overline{C}_{d11}-\overline{C}_{d12}(
\overline {C}_{d22})^-\overline{C}_{d21},
\\
c_{dij}&=&\operatorname{ tr}(B_t\overline{C}_{dij}B_t)
=\operatorname{ tr}(B_tC_{dij}B_t),\qquad  1\leq i,j\leq2.
\end{eqnarray*}
Define $\hat{C}_{dij}=\sum^n_{u=1}\hat{C}_{uij}$, where $\hat
{C}_{uij}=G_i'AG_j$ with $G_1=\hat{T}_u$ and $G_2=\hat{F}_u$. Since
$B\geq0$, we have
%
\begin{eqnarray}\label{eqn2032}
(B_tC_{dij}B_t )_{1\leq i,j\leq2} &=& \pmatrix{ \hat{C}_{d11} & \hat{C}_{d12}
\vspace*{2pt}\cr
\hat{C}_{d21} & \hat{C}_{d22}}
 - n
 \pmatrix{
\hat{\overline{T}}{}^{\prime}_d
\vspace*{2pt}\cr
\hat{\overline{F}}{}^{\prime}_d}
 B \pmatrix{
\hat{\overline{T}}_d & \hat{\overline{F}}_d}
\nonumber
\\[-8pt]
\\[-8pt]
\nonumber
&\leq& (\hat{C}_{dij} )_{1\leq i,j\leq2}.
\end{eqnarray}
Define $q_{dij}=\operatorname{ tr}(\hat{C}_{dij})$ and $q_{uij}=\operatorname{ tr}(\hat{C}_{uij})$.
Then we have $q_{dij}=\sum^n_{u=1}q_{uij}$. It is easy to see that
$q_{u22}>0$ and hence $q_{d22}>0$, which allow us to define
\[
q_d^* =q_{d11}-\frac{q_{d12}^2}{q_{d22}}.
\]
By (\ref{eqn2032}) we have
\[
\pmatrix{c_{d11} & c_{d12}
\vspace*{2pt}\cr
c_{d21} & c_{d22}
}
\leq \pmatrix{
q_{d11} & q_{d12}
\vspace*{2pt}\cr
q_{d21} & q_{d22}}
,
\]
and then by Lemma~\ref{lemma2122} we have
%
\begin{equation}
\label{eqn2036} c_{d11}-\frac{c_{d12}^2}{c_{d22}}\mathbb{I}_{[c_{d22}>0]}
\leq q_d^*,
\end{equation}
with the equality holds when $\hat{\overline{T}}_d=\hat{\overline
{F}}_d=0$. The latter is achieved by designs which are uniform on
periods. To introduce the following theorem, we define
%
\begin{equation}
\label{eqn2037} y^* =\frac{1}{n}\max_dq_d^*.
\end{equation}

\begin{theorem}\label{thm214}
If $C_d=ny^*B_t/(t-1)$ with $y^*$ defined in (\ref{eqn2037}), then
the design $d$ is $\phi_1$-universally optimal.
\end{theorem}
\begin{pf} By conditions (C.1) and (C.2) of $\Phi$ we have
%
\begin{equation}
\label{eqn2033} \Phi(C_d)\leq \Phi(\overline{C}_d),
\end{equation}
where the equality holds if $C_d$ is completely symmetric, that is,
$C_d=\operatorname{ tr}(C_d)B_t/ (t-1)$ since $C_d$ has row and column sums as zero. The
theorem is proved in view of (\ref{eqn2034}), (\ref{eqn2035}),
(\ref{eqn2036}), (\ref{eqn2037}), (\ref{eqn2033}) and condition (C.3).
\end{pf}

\subsection{Existence and equivalence}\label{sec33}
Theorem~\ref{thm214} provides a sufficient condition for a design to
be $\phi_1$-universally optimal. A natural question is the following:
does there exist such a design? This section gives a positive answer as
well as its corresponding implications.

\begin{theorem}\label{thm204}
For any symmetric design, we have:
\begin{longlist}[(1)]
\item[(1)] $C_d$ is completely symmetric;
\item[(2)] $\operatorname{ tr}(C_d)=q_d^*$;
\item[(3)] given any design $d$ there always exist a corresponding
symmetric design which has the same value of $q_d^*$.
\end{longlist}
\end{theorem}

\begin{remark}
Note that Theorem~\ref{thm204} does not hold if we replace $C_d$
therein by $C_d(\tau,l)$. Hence the argument cannot be applied to
$\Phi(C_d(\tau,l))$ directly. This is why we work on $\phi_1$
instead of $\phi_0$ directly.
\end{remark}

\begin{corol}\label{cor204}
\textup{(i)} There exists a symmetric $\phi_1$-universally optimal design $d$ with
%
\begin{equation}
\label{eqn204} C_d=\frac{ny^*B_t}{t-1}.\vspace*{-6pt}
\end{equation}
\begin{longlist}
\item[(ii)] If a design $d$ is $\phi_1$-universally optimal (or $\phi_1$-optimal with $\Phi$ strictly concave or increasing), then we have
(\ref{eqn204}).
\end{longlist}
\end{corol}
\begin{pf}
(i) is proved by Theorems~\ref{thm214} and \ref
{thm204}. (ii) is proved by (i) and the remark in Kiefer's (\citeyear{Kie75})
Proposition 1.
\end{pf}

\subsection{\texorpdfstring{A necessary condition for $\phi_1$-universal optimality}
{A necessary condition for phi1-universal optimality}}\label{sec34}

In this section we give a necessary condition for a design to be $\phi_1$-universally optimal and define quantities that will be useful for
presenting the necessary and sufficient conditions for $\phi_1$-universal optimality in Section~\ref{sec5}. Now define the
function $q_s(x)=q_{s11}+2q_{s12}x+q_{s22}x^2$ and
$q_d(x)=q_{d11}+2q_{d12}x+q_{d22}x^2$. Since $q_{dij}=\sum^n_{u=1}q_{uij}=n\sum_{s \in{\cal S}}p_sq_{sij}$ we have
%
\begin{equation}
\label{eqn2163} q_d(x)=n\sum_{s \in{\cal S}}p_sq_s(x).
\end{equation}
Since $q_{d22}>0$, by direct calculation we have
%
\begin{eqnarray}\label{eqn20311}
q_d^*&=&\min_xq_d(x)
\nonumber
\\[-8pt]
\\[-8pt]
\nonumber
&=&n\min_x\sum_{s \in{\cal S}}p_sq_s(x).
\end{eqnarray}
By (\ref{eqn2037}) and (\ref{eqn20311}) we have
%
\begin{equation}
\label{eqn2164} y^*= \max_P \min_x\sum
_{s\in{\cal S}}p_sq_s(x).\vadjust{\goodbreak}
\end{equation}

Let $d^*$ be a design which maximizes $q^*_d$. By (\ref{eqn2163}),
(\ref{eqn20311}) and (\ref{eqn2164}) we have $\min_x
q_{d^*}(x)=q^*_{d^*}=ny^*$. Since $q_{d22}>0$ the equation
$q_{d^*}(x)=ny^*$ has a unique solution which is denoted by $x^*$. Define
\[
{\cal T}=\bigl\{s\in{\cal S}\dvtx y^*=q_s\bigl(x^*\bigr)\bigr\}.
\]
Lemma~\ref{lemma215} shows that any universally optimal design is
supported on ${\cal T}$.
%
\begin{lemma}\label{lemma215}
If a design $d$ is $\phi_1$-universally optimal (or $\phi_1$-optimal
with $\Phi$ strictly concave or increasing) then we have
\[
p_s=0,\qquad s\notin{\cal T}.
\]
\end{lemma}
\begin{pf} By Corollary~\ref{cor204}, we have $\operatorname{ tr}(C_d)=ny^*$ and
$C_d=\overline{C}_d$. By (\ref{eqn2034}), (\ref{eqn2035}) and
(\ref{eqn2036}) we have $\operatorname{ tr}(\overline{C}_d)\leq q^*_d$. The theorem
is proved in view of (\ref{eqn2037}) and Section~4.4 of \citet{Kus97N2}.
\end{pf}

\subsection{Determination of $x^*$, $y^*$ and ${\cal T}$}\label{sec35}
For a sequence $s=(t_1,t_2,\ldots,t_p)$, define $s_k=(t_1,\ldots,t_k)$ to be
the first $k$ periods of $s$. Particularly, we have $s=s_p$. For $1\leq
k\leq p$ and $1\leq i\leq t$, we define the treatment/sequence index
$f_{s_k,i}=\sum^k_{j=1}1_{t_j=i}$.
To introduce the following theorem, we define two special symmetric
blocks. The symmetry block $\langle di\rangle$ consists of all
sequences having distinct treatments in the $p$ periods. The symmetry
block $\langle re\rangle$ consists of all sequences having distinct
treatments in the first $p - 1$ periods, with the treatment in period
$p - 1$ repeating in period $p$.

\begin{theorem}\label{thm2042}
For any integer $k>t$, define $z_k$ and $r_k$ to be integers satisfying
$k=z_kt+r_k$ and $0< r_k\leq t$.
\begin{longlist}[(iii)]
\item[(i)] If $m>t$ and
%
\begin{equation}
\label{eqn217} \sum^p_{k=m}
\alpha_k \bigl[k\bigl(mt-t^2+1-k\bigr)+t-r_k(t-r_k+1)
\bigr]\geq 0,
\end{equation}
then
\begin{eqnarray*}
x^*&=&0,
\\
y^*&=&\sum^p_{k=m}\alpha_k
\bigl[k(1-1/t)-r_k(t-r_k)/pt\bigr],
\\
{\cal T}&=&\{s\dvtx f_{s_k,i}=z_k\mbox{ or
}z_k+1, 1\leq i\leq t, m\leq k\leq p\}.
\end{eqnarray*}
\item[(ii)] If $p\leq t$ and
%
\begin{equation}
\label{eqn2045} \sum^{p-1}_{k=m}
\alpha_k(k-1) (p+1/t-k)\leq \alpha_p
\bigl[(p-1)^2-(1+1/t)p+1/t\bigr],
\end{equation}
then
\begin{eqnarray*}
x^*&=&1/(p-1),
\\
y^*&=&\sum^p_{k=m}\alpha_k(k-1)
\biggl(1-\frac
{2p-1-k+1/t}{k(p-1)^2} \biggr),
\\
{\cal T}&=&\langle re\rangle\cup\langle di\rangle.
\end{eqnarray*}
When the two sides of (\ref{eqn2045}) are equal, we have ${\cal
T}=\langle re\rangle$.

\item[(iii)] Let
\[
x_0=\frac{\sum^{p-1}_{k=m}\alpha_k(k-1)}{\sum^p_{k=m}\alpha_k(k-1)(k-1-1/t)}.
\]
If $(p-1)^{-1}< x_0<(p-2)^{-1}$, then
\begin{eqnarray*}
x^*&=&x_0,
\\
y^*&=&\sum^p_{k=m}\alpha_k(k-1)
(1-1/k-1/kt)x^2_0-2\sum^{p-1}_{k=m}
\alpha_k(1-1/k)x_0
\\
&&{}+\sum^p_{k=m}\alpha_k(k-1)-2/p,
\\
{\cal T}&=&\langle re\rangle.
\end{eqnarray*}
\end{longlist}
\end{theorem}

\begin{remark}
Under complete experiment, Theorem~\ref{thm2042}(i) applies to the
case $p>t$, and Theorem~\ref{thm2042}(ii) applies to the case
$p\leq t$. Actually Theorem~\ref{thm2042}(i), (ii) reduce to
Theorem 1 of \citet{Kus98}. One can extrapolate by continuity that
Theorem~\ref{thm2042}(i), (ii) cover the cases when the dropout
issue is not very serious.
\end{remark}

\begin{remark}\label{rem416}
When $m=p-1$, we would also discuss parts (i) and (ii) of Theorem
\ref{thm2042}. For (i), a sufficient condition for (\ref{eqn217})
is $p> t+3$. For (ii), inequality~(\ref{eqn2045}) simplifies to
%
\begin{equation}
\label{eqn216}
\alpha_p\geq\frac{(p-2)(1+1/t)}{(p-1)^2-2-1/t}.
\end{equation}
The right-hand side of (\ref{eqn216}) mainly depends on $p$, and it
will become very small for large $p$. Particularly, a sufficient
condition for (\ref{eqn216}) is
\[
a_p\geq\frac{n}{n+1}\frac{(p-2)(1+1/t)}{(p-1)^2-2-1/t}+\frac
{1}{n+1}.
\]
\end{remark}

\section{\texorpdfstring{Linear equations for $\phi_1$-universal optimality}
{Linear equations for phi1-universal optimality}}\label{sec5}
Built upon the results of Section~\ref{sec3}, this section provides
feasible equivalent conditions in approximate design theory for $\phi_1$-universal optimality.\vadjust{\goodbreak}

\subsection{Equations for general designs}\label{sec51}

Recall that $\hat{T}_u=T_uB_t$ and $\hat{F}_u=F_uB_t$, and then we define
%
\begin{eqnarray}\label{eqn2124}
\check{C}_d&=&\check{C}_{d11}-\check{C}_{d12}
\check {C}_{d22}^-\check{C}_{d21},
\nonumber
\\[-8pt]
\\[-8pt]
\nonumber
\check{C}_{dij}&=&\sum^n_{u=1}
\check{C}_{uij},\qquad 1\leq i,j\leq2,
\end{eqnarray}
where
\begin{eqnarray*}
\check{C}_{u11}&=&T_u'(A-B)T_u+
\hat{T}_u'B\hat{T}_u, \qquad \check
{C}_{u12}=T_u'(A-B)F_u+
\hat{T}_u'B\hat{F}_u,
\\
\check{C}_{u21}&=&\check{C}_{u12}', \qquad \check
{C}_{u22}=F_u'(A-B)F_u+
\hat{F}_u'B\hat{F}_u.
\end{eqnarray*}
We shall replace $\check{C}_{uij}$ with $\check{C}_{sij}$ in
emphasizing sequence $s$ instead of subject $u$ of a design. By direct
calculation we have
%
\begin{equation}
\label{eqn212} \check{C}_{dij}=C_{dij}+nG_i'BG_j,\qquad
1\leq i,j\leq2,
\end{equation}
where $G_1=\hat{\overline{T}}_d$ and $G_2=\hat{\overline{F}}_d$.
The following lemma is crucial for the proof of Theorem~\ref{thm2142}.
%
\begin{lemma}\label{lemma212}
If $d$ is $\phi_1$-universally\vspace*{1pt} optimal (or $\phi_1$-optimal with
$\Phi$ strictly concave or increasing), we have
$C_d=\check{C}_d=ny^*B_t/(t-1)$.
\end{lemma}
\begin{pf}
By (\ref{eqn212}) and Lemma~\ref{lemma2122} we have
%
\begin{equation}
\label{eqn2125} C_d\leq \check{C}_d.
\end{equation}
By Corollary~\ref{cor204}(ii) we have
%
\begin{equation}
\label{eqn2123} C_d=ny^*B_t/(t-1).
\end{equation}
Let $\bar{d}$ be the symmetrized version of design $d$ as defined by
(\ref{eqn2122}), and then by (\ref{eqn2124}) we have
%
\begin{equation}
\label{eqn2126} \sum_{\sigma\in{\cal P}}S_{\sigma}'
\check{C}_{dij}S_{\sigma
}/|{\cal P}|=\check{C}_{\bar{d}ij}.
\end{equation}
Again by (\ref{eqn212}) we have $\check{C}_{dij}=G_i'\Lambda G_j$
with $G_1=(\hat{\overline{T}}{}^{\prime}_d,T_d')'$, $G_2=(\hat{\overline
{F}}{}^{\prime}_d,F_d')'$, and
\[
\Lambda=\pmatrix{ nB & 0
\vspace*{2pt}\cr
0 & V}.
\]
Since $\Lambda\geq0$ we have by Proposition 1 of \citet{KunMar00} that
%
\begin{equation}
\label{eqn2129} \sum_{\sigma\in{\cal P}}S_{\sigma}'
\check{C}_dS_{\sigma}/|{\cal P}|\leq\check{C}_{\bar{d}},
\end{equation}
in view of (\ref{eqn2126}). Since $\hat{\overline{T}}_{\bar
{d}}=\hat{\overline{F}}_{\bar{d}}=0$ for the symmetric design $\bar
{d}$, we have
%
\begin{equation}
\label{eqn2128} \check{C}_{\bar{d}}=C_{\bar{d}}.
\end{equation}
Combining (\ref{eqn2125})--(\ref{eqn2129}) and (\ref{eqn2128}), we have
\begin{eqnarray*}
\frac{ny^*}{t-1}B_t=C_d&=&\overline{C}_d
\\
&\leq& \sum_{\sigma\in{\cal P}}S_{\sigma}'
\check{C}_dS_{\sigma
}/|{\cal P}|
\\
&\leq&C_{\bar{d}}.
\end{eqnarray*}
Hence we have $C_{\bar{d}}=ny^*B_t/(t-1)$ in view of Corollary \ref
{cor204} and thus
\[
\sum_{\sigma\in{\cal P}}S_{\sigma}'
\check{C}_dS_{\sigma}/|{\cal P}|= ny^*B_t/(t-1),
\]
which in turn yields
%
\begin{equation}
\label{eqn2127} \operatorname{ tr}(\check{C}_d)=ny^*.
\end{equation}
The lemma is now proved in view of (\ref{eqn2125}) and (\ref{eqn2127}).
\end{pf}

\begin{theorem}\label{thm2142}
A design $d$ is $\phi_1$-universally optimal (or $\phi_1$-optimal
with $\Phi$ strictly concave or increasing) if and only if
%
\begin{eqnarray}
\sum_{s\in{\cal T}}p_s\bigl[
\check{C}_{s11}+x^*\check {C}_{s12}B_t\bigr]&=&
\frac{y^*}{t-1}B_t,\label{eqn11}
\\
\sum_{s\in{\cal T}}p_s\bigl[
\check{C}_{s21}+x^*\check {C}_{s22}B_t\bigr]&=&0,
\label{eqn22}
\\
\sum_{s\in{\cal T}}p_sB\bigl(
\hat{T}_s+x^*\hat{F}_s\bigr)&=&0,\label {eqn33}
\\
\sum_{s\in{\cal T}}p_s&=&1,
\\
p_s&=& 0,\qquad s\notin{\cal T}.\label{eqn44}
\end{eqnarray}
\end{theorem}

Based on Theorem~\ref{thm214} and Corollary~\ref{cor204}, (\ref
{eqn204}) is also a necessary and sufficient condition for $\phi_1$-universal optimality.
However, (\ref{eqn204}) is not directly
applicable for identifying designs. Note that the conditions in Theorem
\ref{thm2142} are merely linear equation systems for $p_s$, and hence
can be easily implemented to derive exact designs. See Section~\ref{sec6}.

\subsection{Equations for symmetric designs}

Note that $q_s(x)$ is invariant to treatment permutation, that is,
%
\begin{equation}
\label{eqn2042} q_s(x)= q_{\sigma s}(x).\vadjust{\goodbreak}
\end{equation}
Combining Theorem 4.5 of \citet{Kus97N2}, Theorem~\ref{thm204},
Corollary~\ref{cor204}, Lemma~\ref{lemma215} and equation (\ref
{eqn2042}), we have the following:

\begin{theorem}\label{thm416}
A symmetric design is $\phi_1$-universally optimal if
\begin{eqnarray*}
\sum_{s\langle\in\rangle{\cal T}}p_{\langle s\rangle}q_s'
\bigl(x^*\bigr)&=&0,
\\
\sum_{s\langle\in\rangle{\cal T}}p_{\langle s\rangle}&=&1,
\\
p_s&=&0, \qquad s\notin{\cal T},
\end{eqnarray*}
where $q_s'(x)$ is the derivative of $q_s(x)$ with respective to $x$.
\end{theorem}

\section{Exact designs}\label{sec6}

This section gives algorithms to identify efficient exact designs based
on the optimality equations in Section~\ref{sec5}. Results are
compared to designs proposed in literature. For the matrix $C_d(\tau
,l)$, denote its eigenvalues by $0=\lambda_1\leq\lambda_2\leq
\cdots\leq\lambda_t$. We define the criteria of $A$, $D$, $E$ and $T$ as:
\begin{itemize}
\item$\Phi_A(C_d(\tau,l))=(t-1)(n\sum^t_{i=2}\lambda_i^{-1})^{-1}$. [$\lambda_2=0$ implies $\Phi_A(C_d(\tau,l))=0$];\vspace*{1pt}
\item$\Phi_D(C_d(\tau,l))=n^{-1}(\prod^t_{i=2}\lambda_i)^{1/(t-1)}$;
\item$\Phi_E(C_d(\tau,l))=n^{-1}\lambda_2$;
\item$\Phi_T(C_d(\tau,l))=[n(t-1)]^{-1}(\sum^t_{i=2}\lambda_i)$.
\end{itemize}
Section~\ref{sec61} provides an algorithm to derive exact designs for
general configurations of $p,t,n$. Section~\ref{sec62} illustrates
how to derive symmetric designs by straightforward calculations. In
utilizing Lemma~\ref{lemma203}, $e_1(d)$ is further bounded by
$\tilde{e}_1=\phi_1(d)/\phi_1(\tilde{d})$, where $\tilde{d}$ is a
$\phi_1$-optimal design in asymptotic design theory which may not
necessarily exist as an exact design. Thus the function $\ell
(d)=\tilde{e}_1(d)g(d)$ serves as a feasible lower bound of $e_0(d)$.

\subsection{General exact designs}\label{sec61}
This section gives an algorithm to derive efficient exact optimal
designs for any given configuration of $p,t,n$ and compares them to
designs in literature. Note that the latter designs are proposed for
judiciously chosen $p,t,n$ while our algorithm works for any
configuration of $p,t,n$. Even under these chosen circumstances our
designs are still shown to be more efficient and robust. By Theorem
\ref{thm2142} we have the following:
%
\begin{corol}
A design $d$ is $\phi_1$-universally optimal (or $\phi_1$-optimal
with $\Phi$ strictly concave or increasing) if and only if
%
\begin{eqnarray}
\sum_{s\in{\cal T}}n_s\bigl[
\check{C}_{s11}+x^*\check {C}_{s12}B_t\bigr]&=&
\frac{ny^*}{t-1}B_t,\label{eqn327}
\\
\sum_{s\in{\cal T}}n_s\bigl[
\check{C}_{s21}+x^*\check {C}_{s22}B_t\bigr]&=&0,
\label{eqn3272}
\\
\sum_{s\in{\cal T}}n_sB\bigl(
\hat{T}_s+x^*\hat{F}_s\bigr)&=&0,\label {eqn3273}
\\
\sum_{s\in{\cal T}}n_s&=& n,\label{eqn3274}
\\
n_s&=&0, \qquad s\notin{\cal T}.\label{eqn3275}
\end{eqnarray}
\end{corol}
Note that an exact design satisfying equations (\ref{eqn327})--(\ref
{eqn3275}) does not necessarily exist due to the discrete nature of
the problem, especially when the dropout mechanism is arbitrary.
However, as shown by the following examples, it is plausible to find a
design which is as \textit{close} to satisfying equations (\ref
{eqn327})--(\ref{eqn3275}) as possible. Specifically, let $N_{\cal
T}=\{n_s,s\in{\cal T}\}'$, and then equations (\ref{eqn327})--(\ref
{eqn3273}) could be written in a matrix form as
\[
\label{eqn3276} X_{\cal T}N_{\cal T}=Y_{\cal T},
\]
with $X_{\cal T}$ and $Y_{\cal T}$ uniquely determined by equations
(\ref{eqn327})--(\ref{eqn3275}) and the ordering of the $n_s$ in
the vector $N_{\cal T}$. To find an efficient design for an arbitrarily
given $n$, one could choose a design which

Minimizes
%
\begin{equation}
\label{eqn329}\Vert X_{\cal T}N_{\cal T}-Y_{\cal T}\Vert ,
\end{equation}

subject to
\[
1_{|{\cal T}|}'N_{\cal T}=n.
\]
Here $\Vert \cdot\Vert $ is a norm for a vector. For all subsequent examples
in this section, we take $\Vert \cdot\Vert $ to be the Euclidean norm. Then
solving for (\ref{eqn329}) is straightforward by utilizing integer
optimization packages/softwares. Note that the computational complexity
of the above minimization problem depends on $|{\cal T}|$, which in
turn depends on $p$ and $t$ only.

Besides maximizing the expectation $\phi_0(d)=\E\Phi(C_d(\tau,l))$,
one might also be interested in minimizing the variance $V_{\Phi
}(d)=\operatorname{Var}(\Phi(C_d(\tau,l)))$ to achieve robustness. To compare two
designs under these two functions, we define $\phi_0(d,d')=\phi_0(d)/\phi_0(d')$ and $V_{\Phi}(d,d')=V_{\Phi}(d)/V_{\Phi}(d')$.

\subsubsection{\texorpdfstring{Comparisons to designs of Low, Lewis and Prescott (\citeyear{LowLewPre99})}
{Comparisons to designs of Low, Lewis and Prescott (1999)}}
The setup and target of \citet{LowLewPre99} are the same as
in this paper. However, they searched all combinations of Latin squares
for the special cases of $p=t=4,n=16$ and $p=t=4,n=24$ only.

%
\begin{figure}

\includegraphics{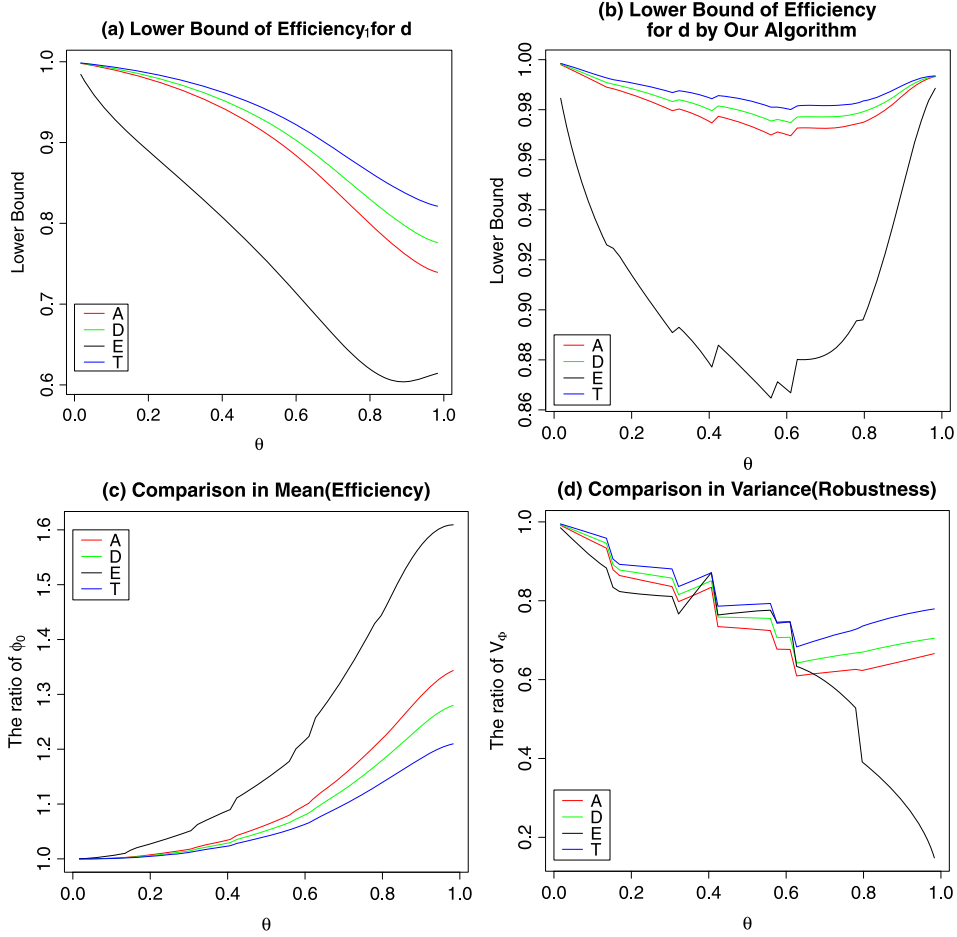}

\caption{The letters $A$, $D$, $E$ and $T$ represent the choice of criteria
function $\Phi$.
\textup{(a)} The lower bound of efficiency $\ell(d_1)$ for
$a=(0,0,\theta,1-\theta)$ with $\theta\in(0,1)$.
\textup{(b)} The lower bound of efficiency $\ell(d)$ with $d$
obtained by algorithm (\protect\ref{eqn329}).
Particularly $\theta=1/2$ implies $d=d_2$. \textup{(c)} The ratio of
mean: $\phi_0(d)/\phi_0(d_1)$.
\textup{(d)} The ratio of variance: $V_{\Phi}(d)/V_{\Phi
}(d_1)$.}\label{fig1}
\end{figure}

When $p=t=4$ and $n=16$, they proposed a design as shown by Figure~\ref{fig1}(b)
therein, which is said to be $d_1$ here. By algorithm (\ref{eqn329}),
the dropout mechanism $\vec{a}=(0,0,1/2,1/2)$ yields $d_2$.
\begin{eqnarray*}
d_2\dvtx\qquad  %
\begin{array} {cccccccccccccccc} 2 & 1 & 2 & 3 &
3 & 4 & 3 & 2 & 1 & 1 & 1& 2 & 4 & 4 & 4 & 3
\\
4 & 4 & 3 & 4 & 1 & 1 & 2 & 1 & 2 & 2 & 3 & 4 & 3 & 3 & 2 & 1
\\
3 & 2 & 1 & 1 & 2 & 3 & 4 & 4 & 3 & 3 & 4 & 1 & 2 & 2 & 1 & 4
\\
3 & 2 & 1 & 1 & 2 & 3 & 4 & 4 & 4 & 4 & 2 & 3 & 1 & 1 & 3 & 2
\end{array}
.
\end{eqnarray*}
%

\begin{table}
\caption{Performance of $d_1$ under $\vec{a}=(0,0,1/2,1/2)$}\label{tab3292}
\begin{tabular*}{\textwidth}{@{\extracolsep{\fill}}lccccc@{}}
\hline
$\bolds{\Phi}$ & $\bolds{\phi_0(d_1)}$ & $\bolds{V_{\Phi}(d_1)}$ & $\bolds{\tilde{e}_1(d_1)}$ &
$\bolds{g(d_1)}$ & $\bolds{\ell(d_1)}$\\
\hline
$A$ &0.6646345& 0.07223834& 0.9558432& 0.9592851& 0.9169261\\
$D$ &0.6747419& 0.06776632 &0.9603310& 0.9693223& 0.9308702\\
$E$ &0.5528575 &0.09039916 &0.8960473& 0.8512042& 0.7627192\\
$T$ &0.6848634 &0.06334558 &0.9650531& 0.9790485& 0.9448338\\
\hline
\end{tabular*} \vspace*{-2pt}
\end{table}

%
\begin{table}[b]\vspace*{-2pt}
\caption{Performance of $d_2$ under $\vec{a}=(0,0,1/2,1/2)$}\label{tab329}
\begin{tabular*}{\textwidth}{@{\extracolsep{\fill}}lccccc@{}}
\hline
$\bolds{\Phi}$ & $\bolds{\phi_0(d_2)}$ & $\bolds{V_{\Phi}(d_2)}$ & $\bolds{\tilde{e}_1(d_2)}$ &
$\bolds{g(d_2)}$ & $\bolds{\ell(d_2)}$\\
\hline
$A$& 0.7058735& 0.05266523& 0.9989759 &0.9748175 &0.9738192\\
$D$& 0.7094851 &0.05129209& 0.9991830 &0.9796020 &0.9788017\\
$E$ &0.6337475 &0.06979073& 0.9848636 &0.8877519 &0.8743145\\
$T$ &0.7130567 &0.05005383& 0.9993922 &0.9843273 &0.9837291\\
\hline
\end{tabular*}
\end{table}

Tables~\ref{tab3292} and~\ref{tab329} summarize the performances of
designs $d_1$ and $d_2$ under criteria of $A$, $D$, $E$ and $T$. Since
$e_0(d)\geq\ell(d)=\tilde{e}_1(d)g(d)$, a design $d$ would be $\phi_0$-efficient
if both $\tilde{e}_1(d)$ and $g(d)$ are close to unity.
Algorithm (\ref{eqn329}) focuses on $e_1(d)$ and provides a
satisfactory solution in view of the column of $\tilde{e}_1$ in Table
\ref{tab329}. We observe that the values of $g$ in both of these
tables are very close to unity except for $E$-criterion. Notice that the
values of gap function $g$ for $T$-criterion are always the largest among
all criteria, which is due to the linearity of $T$-criterion.

In comparison, $d_2$ is more efficient and robust than $d_1$ under all
criteria in view of the columns of $\phi_0$ and $V_{\Phi}$,
respectively. A lesson from the latter is that a design with a more
diverse composition of sequences is generally more robust. Here in
$d_2$, only the sequences of $1234$ and $4321$ appear twice while each
of the remaining sequences appears only once. \citet{LowLewPre99} had similar observations.

We now consider the performance of a design obtained by algorithm (\ref
{eqn329}) for dropout mechanisms of the form $\vec{a}=(0,0,\theta
,1-\theta)$, $0<\theta<1$. By heuristic arguments in Section \ref
{sec22}, the value of gap function $g$ would be smaller if there is
larger variability in $l$. This is supported by the $U$-shape curve of
$\ell(d)$ in Figure~\ref{fig1}(b). From Figure~\ref{fig1}(a), we see that the efficiency
of $d_1$ has a reverse relationship with the value of $\theta$. Figure
\ref{fig1}(c) shows that the advantage of our algorithm against $d_1$ is more
obvious when there is large chance of dropout. This means that our
algorithm succeeded in adapting the choice of designs to different
dropout mechanisms. Figure~\ref{fig1}(d) shows that the design by our algorithm
is also more robust than $d_1$ against the randomness of subject dropout.
When $p=t=4$ and $n=24$, \citet{LowLewPre99} proposed a
design which consists of two copies of three distinct $4\times4$ Latin
squares, which is denoted by $d_3$ here. When\vadjust{\goodbreak} $\vec
{a}=(0,1/10,2/5,1/2)$ our algorithm yields $d_4$ which consists of one
copy of the first twelve sequences and two copies of the last six
sequences of (\ref{eqn43}). According to the last two columns of
Table~\ref{tab403}, $d_4$ outperforms $d_3$ in terms of both
efficiency and robustness with the exception for the robustness under
$E$-criterion.
%
\begin{eqnarray}
\label{eqn43} d_4\dvtx\qquad %
\begin{array} {cccccccccccc@{\hspace*{4pt}}|@{\hspace*{4pt}}cccccc}
2 & 2 & 2& 3 & 4& 4 &2& 3 & 3& 4 & 4& 3 & 3& 4 & 2 & 1 & 1 & 1
\\
3 & 4 & 3 & 4 & 3 & 2 & 3 & 4 & 4 & 2 & 2 & 2 & 1 & 1 & 1 & 2 & 4 & 3
\\
4 & 3 & 1 & 1 & 2 & 1 & 1 & 1 & 2 & 1 & 3 & 4 & 2 & 3 & 4 & 4 & 3 & 2
\\
4 & 3 & 1 & 1& 2 & 1& 4 & 2& 1 & 3& 1 & 1& 2 & 3 & 4 & 3 & 2 & 4 \end{array}
\times2.
\end{eqnarray}
%

\begin{table}
\caption{Performance of $d_3$ and $d_4$ under $\vec{a}=(0,1/10,2/5,1/2)$}\label{tab403}
\begin{tabular*}{\textwidth}{@{\extracolsep{\fill}}lccccccc@{}}
\hline
$\bolds{\Phi}$ & $\bolds{\phi_0(d_4)}$ & $\bolds{V_{\Phi}(d_4)}$ & $\bolds{\tilde{e}_1(d_4)}$ &
$\bolds{g(d_4)}$ & $\bolds{\ell(d_4)}$& $\bolds{\phi_0(d_4,d_3)}$ & $\bolds{V_{\Phi}(d_4,d_3)}$\\
\hline
$A$ &0.6791& 0.0526 &0.999983 &0.9777& 0.9777& 1.0112 &0.9705\\
$D$& 0.6822& 0.0516& 0.999983 &0.9822& 0.9821 &1.0115 &0.9631\\
$E$& 0.6118& 0.0648 &0.999979 &0.8809& 0.8809 &1.0089 &1.0386\\
$T$ &0.6852& 0.0506 &0.999983 &0.9866& 0.9869 &1.0118 &0.9562\\
\hline
\end{tabular*}
\end{table}

\subsubsection{\texorpdfstring{Comparison to designs of Bose and Bagchi (\citeyear{BosBag08}),
Majumdar, Dean and Lewis (\citeyear{MajDeaLew08}) and Zhao and Majumdar (\citeyear{ZM12})}
{Comparison to designs of Bose and Bagchi (2008),
Majumdar, Dean and Lewis (2008) and Zhao and Majumdar (2012)}}\label{sec41}

When the realization of subject dropout $l$ is not random, we have
$\phi_0\equiv\phi_1$. In this case, \citet{BosBag08} have the
following results:
\begin{longlist}[(1)]
\item[(1)] When $p=t\geq3$ is a prime or primer power and $n=t(t-1)$,
a design is found to be universally optimal whenever $a_q=1$ for any
$3\leq q\leq p$.
\item[(2)] When $p=t\geq3$ is a prime or primer power, $t\equiv
3$(mod $4$) and $n=2t$, a design is found to be universally optimal
whenever $a_q=1$ with $q=(p+1)/2$ or~$p$.
\item[(3)] When $p=t\geq3$ is a prime or primer power, $t\equiv
1$(mod $4$) and $n=4t$, a~design is found to be universally optimal
whenever $a_q=1$ with $q=(p+1)/2$ or~$p$.
\end{longlist}

For example, when $t=p=5$ the smallest $n$ should be $4t=20$. In this
case the design proposed by them is universally optimal, either when
the experiment is complete or when all subjects immediately drop out
after period $3$ with probability~$1$, that is, $a_3=1$. We denote this
design by $d_5$ which is given by Example~3 of \citet{BosBag08}.
When $\vec{a}=(0,1/20,3/20,1/5,3/5)$ algorithm (\ref{eqn329}) yields
$d_6$ as follows:
\begin{eqnarray*}
d_6\dvtx\qquad
\begin{array} {cccccccccccccccccccc} 1 & 2 & 4 &
4 & 3 & 2 & 1 & 2 & 1 & 1 & 3 & 2 & 4 & 5 & 5 & 5 & 5 & 3 & 4 & 3
\\
2 & 5 & 1 & 2 & 1 & 3 & 2 & 3 & 5 & 4 & 4 & 5 & 5 & 4 & 4 & 2 & 3 & 1 & 3 & 1
\\
3 & 4 & 3 & 5 & 5 & 4 & 4 & 5 & 3 & 3 & 1 & 1 & 1 & 2 & 2 & 1 & 2 & 2 & 5 & 4
\\
4 & 1 & 5 & 3 & 2 & 5 & 3 & 1 & 2 & 2 & 5 & 3 & 3 & 1 & 1 & 4 & 4 & 4 & 2 & 5
\\
4 & 1 & 5 & 3 & 2 & 1 & 5 & 4 & 4 & 5 & 2 & 4 & 2 & 3 & 3 & 3 & 1 & 5 & 1 & 2
\end{array} %
.
\end{eqnarray*}
Table~\ref{tab416} shows that $d_6$ is more efficient and robust than
$d_5$ under criteria of $A$, $D$ and $T$, while the result is reversed under
the criterion of $E$. The reason for the latter is that $d_5$ did a
better job in avoiding disconnected designs under subject dropout, that
is, $\Phi_E(C_d(\tau,l))=0$.

\begin{table}
\caption{Performance of $d_5$ and $d_6$ under $\vec
{a}=(0,1/20,3/20,1/5,3/5)$}\label{tab416}
\begin{tabular*}{\textwidth}{@{\extracolsep{\fill}}lccccccc@{}}
\hline
$\bolds{\Phi}$ & $\bolds{\phi_0(d_6)}$ & $\bolds{V_{\Phi}(d_6)}$ & $\bolds{\tilde{e}_1(d_6)}$ &
$\bolds{g(d_6)}$ & $\bolds{\ell(d_6)}$& $\bolds{\phi_0(d_6,d_5)}$ & $\bolds{V_{\Phi}(d_6,d_5)}$\\
\hline
$A$& 0.7555& 0.05944 &0.99888& 0.97848& 0.97738 &1.00117& 0.98842\\
$D$ &0.7589& 0.05827 &0.99891& 0.98277& 0.98170& 1.00172& 0.98333\\
$E$& 0.6712& 0.07399 &0.99091& 0.87621& 0.86825& 0.99449& 1.03435\\
$T$ &0.7621& 0.05719 &0.99894& 0.98700& 0.98595& 1.00224& 0.97877\\
\hline
\end{tabular*}
\end{table}

\begin{table}[b]
\caption{Performance of $d_7$ and $d_8$ under $\vec
{a}=(0,0,1/3,1/3,1/3)$}\label{tab4162}
\begin{tabular*}{\textwidth}{@{\extracolsep{\fill}}lccccccc@{}}
\hline
$\bolds{\Phi}$ & $\bolds{\phi_0(d_8)}$ & $\bolds{V_{\Phi}(d_8)}$ & $\bolds{\tilde{e}_1(d_8)}$ &
$\bolds{g(d_8)}$ & $\bolds{\ell(d_8)}$& $\bolds{\phi_0(d_8,d_7)}$ & $\bolds{V_{\Phi}(d_8,d_7)}$\\
\hline
$A$ &1.2340 &0.053908 & 1 &0.99591& 0.99591& 1.11018 &0.57705\\
$D$ &1.2347 &0.053736 & 1 &0.99643& 0.99643 &1.10598 &0.59362\\
$E$ &1.2004 &0.059782 & 1 &0.96877& 0.96877 &1.16339 &0.51397\\
$T$& 1.2353& 0.053573 & 1 &0.99696& 0.99696 &1.10177 &0.60992\\
\hline
\end{tabular*}
\end{table}

Since the magnitude of the differences between $d_5$ and $d_6$ are
small in terms of both efficiency and robustness, we conclude that the
designs of \citet{BosBag08} successfully defended the loss of
information due to subject dropout. The same conclusion applies to
Majumdar, Dean and Lewis (\citeyear{MajDeaLew08}) and Zhao and Majumdar (\citeyear{ZM12}) since they
use similar ideas.

\subsubsection{\texorpdfstring{Comparisons to designs of Kushner (\citeyear{Kus98})}{Comparisons to designs of Kushner
(1998)}}
Kushner (\citeyear{Kus98}) derived conditions for universal optimality as a special
case of ours under complete experiment. Particularly, when $t=3$, $p=5$
and $n=30$, Example 4 of \citet{Kus98} gives a design satisfying the
optimality equations therein, which is denoted $d_7$ here. When $\vec
{a}=(0,0,1/3,1/3,1/3)$ our algorithm gives $d_8$ which consist of five
copies of (\ref{eqn416}),

\begin{eqnarray}
\label{eqn416} d_8\dvtx\qquad
\begin{array} {cccccc} 1 & 2 & 3 &
3 & 1 & 2
\\
3 & 3 & 2 & 1 & 2 & 1
\\
2 & 1 & 1 & 2 & 3 & 3
\\
2 & 1 & 1 & 2 & 3 & 3
\\
1 & 2 & 3 & 3 & 2 & 1 \end{array} %
\times5.
\end{eqnarray}

Based on Table~\ref{tab4162} $d_8$ outperforms $d_7$ in terms of both
efficiency and robustness even though $d_7$ is universally optimal
under complete experiment.
%

\subsection{Symmetric exact designs}\label{sec62}
This section illustrates the usage of Theorem~\ref{thm416} in
deriving efficient symmetric exact designs. By Remark~\ref{rem416} in
Section~\ref{sec35}, when $t=2$, $p=6$ and $m=p-1=5$, inequality
(\ref{eqn217}) in Theorem~\ref{thm2042} always holds regardless of
the value of $\vec{a}$. By applying Theorem~\ref{thm2042}(i), we
have $x^*=0$ and hence $q_s'(x^*)=2q_{s12}$. Moreover, it is easy to
see that the support ${\cal T}$ essentially contains all sequences
which assign a subject to each of the two treatments for $3$ out of the
total of $6$ periods, and hence $|{\cal T}|=20$. Within each symmetric
block, there are two sequences since $t=2$. Hence there are $10$
symmetric blocks. However, it is not necessary to include all these
symmetric blocks in the design. Particularly when $\vec
{a}=(0,0,0,0,2/5,3/5)$, we have
$q_{s_1}'(x^*)/q_{s_2}'(x^*)=q_{s_112}/q_{s_212}=-6.01$ for
$s_1=122121$ and $s_2=122211$. In the spirit of Theorem~\ref{thm416}
we propose a small sized design, $d_9$, which consists of one copy of
sequences $122121$ and $211212$ and six copies of the sequences
$122211$ and $211122$. So we have $n=14$ for $d_9$. The point is that
we have the freedom of selecting different subclasses of ${\cal T}$.
The performance of $d_9$ is given in Table~\ref{tab417}. It shows the
high efficiency and robustness of~$d_9$. Note that when $t=2$ all
criteria are equivalent.

\begin{table}
\caption{Performance of $d_9$ under $\vec
{a}=(0,0,0,0,2/5,3/5)$}\label{tab417}
\begin{tabular*}{\textwidth}{@{\extracolsep{\fill}}lccccc@{}}
\hline
$\bolds{\Phi}$ & $\bolds{\phi_0(d_9)}$ & $\bolds{V_{\Phi}(d_9)}$ & $\bolds{\tilde{e}_1(d_9)}$ &
$\bolds{g(d_9)}$ & $\bolds{\ell(d_9)}$\\
\hline
$A$, $D$, $E$, $T$ and etc.& 2.7368& 0.09152 &0.99511 &0.997823 &0.99295\\
\hline
\end{tabular*}
\end{table}

\section{Discussions}\label{sec7}
Subject dropout is a very important issue in planning a crossover
design. It is shown by Table~\ref{tab4162} and other examples in
literature that an optimal design under complete experiment is no
longer optimal and possibly even disconnected when there is subject
dropout. However, the problem has received very limited attention in
literature so far, and the majority of the research assumes that there
is no subject dropout. \citet{BosBag08}, Majumdar, Dean and Lewis
(\citeyear{MajDeaLew08}), Zhao and Majumdar (\citeyear{ZM12}) all considered the nested structure
such that a design, together with its subdesign, obtained by taking
only the first $q$($<p$) periods, are both optimal or efficient.
Naturally such designs would still be efficient when all subjects drop
out at periods between $p$ and $q$. The issue with this approach is
that we lose adaptation to different dropout mechanisms. Furthermore,
their methods only apply to special configurations of $p,t,n$.

In order to take into account the dropout mechanism, one has to make
assumptions to formulate the dropout mechanism. This paper adopts two
mild assumptions and works on the target function $\phi_0$ which is
given by taking the expectation of a regular optimality criterion with
respect to a given dropout mechanism. Actually \citet{LowLewPre99}
have followed the same approach. However,\vadjust{\goodbreak} they only provided two
case studies, and there were no theoretical results regarding how to
identify an efficient design in general. The latter problem is itself
intractable. To tackle it, we propose to use the surrogate target
function of $\phi_1$ in place of $\phi_0$. It turns out that this
replacement is very successful. Examples in Section~\ref{sec6} show
that $\phi_1$-optimal (or highly efficient) designs are also highly
efficient under $\phi_0$. Moreover, these designs are also shown to be
very robust against the randomness of subject dropout due to the
substantial diversity in the composition of treatment sequences.

Theoretically, we derive feasible, equivalent conditions for a design
to be $\phi_1$-universally optimal in asymptotic design theory. These
conditions are essentially linear equations with respect to proportions
of treatment sequences from ${\cal T}$, a subclass of all possible
treatment sequences. A solution for the equations, which yields an
exact design, does not necessarily exist due to the discrete nature of
the problem. However, one can follow the spirit of the conditions and
easily propose an applicable algorithm to derive an efficient exact
design for any criterion and any configuration of $p,t,n$. In this
paper, we adopt algorithm (\ref{eqn329}) for general designs as well
as the approach in Section~\ref{sec62} for symmetric designs.

The problem of identifying exact designs for large values of $p$ and
$t$ remains as an open problem. The critical difficulty is that as $p$
and $t$ grow the size of the support for admissible sequences, $|{\cal
T}|$, increases very fast. Typically ${\cal T}$ contains two distinct
symmetric blocks, in which case $p=t=6$ usually yields $|{\cal
T}|=2\times6!=1440$. That means the majority of the sequences in
${\cal T}$ would not appear in the design for a moderate value of $n$.
The same issue has appeared in \citet{Kus97N2}. If we adopt the
approach of symmetric designs as in Section~\ref{sec62} we would need
$n$ to be as of the same magnitude as $|{\cal T}|$. On the other hand,
algorithm (\ref{eqn329}) is essentially an integer programming
problem and the number of the integer variables is equal to $|{\cal
T}|$. Hence it would be infeasible for a computer to handle when
$|{\cal T}|$ is too large. For this problem, one possible solution is
to reduce the size of ${\cal T}$ through the study of intrinsic
relationships among treatment sequences. Another approach is to resort
to algorithm improvement.\vspace*{-1pt}

\section{Proofs}\label{sec8}\vspace*{-1pt}
\mbox{}
\begin{pf*}{Proof of Lemma~\ref{lemma125}} It would be enough to show that
$V=\E O$. First, it is easy to show that
$B^{m_1}_{ij}B^{m_2}_{jk}=B^{\mathrm{min}(m_1, m_2)}_{ik}$. We have
$MU=\operatorname{diag}(1_{l_1},1_{l_2},\ldots,\break 1_{l_n})$ and
$MZ=({I_{l_1p}^{l_1{}^\prime}},{I_{l_2p}^{l_2{}^\prime}},\ldots, {I_{l_np}^{l_n{}^\prime}})'$.
Then we have
\begin{eqnarray*}
\operatorname{pr}^{\perp}(MU)&=&\operatorname{diag}(B_{l_1},\ldots,B_{l_n}),
\\[-2pt]
\operatorname{pr}^{\perp}(MU)MZ&=&\bigl({B_{l_1p}^{l_1{}^\prime}},{B_{l_2p}^{l_2{}^\prime}},
\ldots, {B_{l_np}^{l_n{}^\prime}}\bigr)',
\\[-2pt]
Z'M'\operatorname{pr}^{\perp}(MU)MZ&=&\sum
^n_{u=1}B_{p}^{l_u}=\sum
^p_{i=1}h_iB_p^i.\vadjust{\goodbreak}
\end{eqnarray*}
Without loss of generality, we could assume $h_p>0$. Then one choice of
the $g$-inverse of $Z'M'\operatorname{pr}^{\perp}(MU)MZ$ is $\sum^p_{i=1}g_iB_p^i$ where
%
\begin{eqnarray}
g_i&=&r_i^{-1}-r_{i+1}^{-1},\qquad 1
\leq i\leq p-1,\label{eqn910}
\\
g_p&=&h_p^{-1},\label{eqn9102}
\end{eqnarray}
with $h_k=\sum^n_{i=1}1_{l_i=k},1\leq k\leq p$, and $r_i=\sum^p_{k=i}h_i$ denotes the number of subjects remaining at period $i$,
$1\leq i\leq p$. Note that if $h_p=0$, the value of $p$ in (\ref
{eqn910}) and~(\ref{eqn9102}) should be replaced by $\tilde{p}=\operatorname{max}\{
k\dvtx h_k>0\}$, and for $k>\tilde{p}$ we let $h_k=0$. It is easily seen
that the following arguments and thus the lemma would still hold. Now
we have
\begin{eqnarray*}
\operatorname{pr}^{\perp}(MZ|MU)&=&\operatorname{pr}^{\perp}(MU)-\operatorname{pr}\bigl(\operatorname{pr}^{\perp}(MU)MZ
\bigr)
\\
&=&\operatorname{diag}(B_{l_1},\ldots,B_{l_n})-\Delta,
\\
\Delta&=& \Biggl(\sum^p_{k=1}g_kB_{l_il_j}^{\min(k,l_i,l_j)}
\Biggr)_{i,j=1,2,\ldots,n}.
\end{eqnarray*}
Let $O=(O_{ij})_{1\leq i,j\leq n}=M'\operatorname{pr}^{\perp}(MZ|MU)M$, and then we have
\begin{eqnarray*}
O_{ii}&=&B_p^{l_i}-\sum
^p_{k=1} g_kB_p^{\mathrm{min}(k,l_i)},
\\
O_{ij}&=&-\sum^p_{k=1}
g_kB_p^{\mathrm{min}(k,l_i,l_j)}.
\end{eqnarray*}
We will derive the expectation of $\sum^p_{k=1} g_kB_p^{\mathrm{min}(k,l_i)}$
and other components could be dealt with by similar arguments. First we
have the decomposition
\begin{eqnarray*}
\sum^p_{k=1} g_kB_p^{\mathrm{min}(k,l_i)}&=&
\sum^{l_i-1}_{k=1}g_kB_p^k+
\sum^p_{k=l_i}g_kB_p^{l_i}
\\
&=&\sum^{l_i-1}_{k=1} \biggl(
\frac{1}{r_k}-\frac{1}{r_{k+1}} \biggr)B_p^k+
\frac{1}{r_{l_i}}B_p^{l_i}.
\end{eqnarray*}
When $k\leq l_i$ and $l_i$ is given, we know that $r_k-1$ follows the
binomial distribution with parameters $n-1$ and $a_{kp}$. Hence we have
\begin{eqnarray*}
\E\bigl(r_k^{-1}|l_i,k\leq l_i
\bigr)&=&\sum^{n-1}_{j=0}\frac{1}{j+1}
\frac
{(n-1)!}{(n-1-j)!j!}a_{kp}^j(1-a_{kp})^{n-1-j}
\\
&=&\frac{1-(1-a_{kp})^n}{na_{kp}}:=b_k.
\end{eqnarray*}
Hence we have
\begin{eqnarray*}
\E \Biggl(\sum^p_{k=1}
g_kB_p^{\mathrm{min}(k,l_i)}\bigg|l_i,1\leq i\leq n
\Biggr)&=&\sum^{l_i-1}_{k=1}(b_k-b_{k+1})B_p^k+b_{l_i}B_p^{l_i}
\\
&=&\sum^p_{k=1} \bigl[(b_k-b_{k+1})1_{[k<l_i]}+b_k1_{[k=l_i]}
\bigr]B_p^k.
\end{eqnarray*}
Here we have the convention of $b_{p+1}=0$ for notational convenience. Hence
\begin{eqnarray*}
\E \Biggl(\sum^p_{k=1}
g_kB_p^{\mathrm{min}(k,l_i)} \Biggr)&=&\sum
^p_{k=1} \bigl[(b_k-b_{k+1})a_{k+1,p}+b_k(a_{kp}-a_{k+1,p})
\bigr]B_p^k
\\
&=&\sum^p_{k=1}(a_{kp}b_k-a_{k+1,p}b_{k+1})B_p^k
\\
&=&\frac{1}{n}\sum^p_{k=1}
\bigl(a_{1k}^n-a_{1,k-1}^n
\bigr)B_p^k.
\end{eqnarray*}
Following this strategy, it is easy to show that
\begin{eqnarray*}
\E O_{ii}&=&\sum^p_{k=1}
\bigl[a_k-n^{-1}\bigl(a_{1k}^n-a_{1,k-1}^n
\bigr) \bigr]B_p^k,
\\
\E O_{ij}&=&-\frac{1}{n}\sum^p_{k=1}
\bigl(a_k+a_{k+1,p}a_{1k}^n-a_{kp}a_{1,k-1}^n
\bigr)B_p^k.
\end{eqnarray*}
Then we have $V=\E O$.
\end{pf*}

\begin{pf*}{Proof of Theorem~\ref{thm204}} By definition of symmetric
designs we have
%
\begin{eqnarray}
T_{\sigma d}&=&T_dS_{\sigma}
\nonumber
\\
&=&(\tilde{S}_{\sigma,d}\otimes I_p)T_d,
\label{eqn2038}
\\
F_{\sigma d}&=&(\tilde{S}_{\sigma,d}\otimes I_p)F_d,
\nonumber
\end{eqnarray}
where $\tilde{S}_{\sigma,d}$ is a permutation matrix for subjects
induced by $\sigma$ and (symmetric)~$d$. Note that we have $(\tilde
{S}_{\sigma,d}\otimes I_p)'V(\tilde{S}_{\sigma,d}\otimes I_p)=V$. So
$C_{dij},1\leq i,j\leq2$, are completely symmetric and hence $C_d$ is
completely symmetric for a symmetric design~$d$. This yields
\[
C_d=\overline{C}_d,
\]
and the equality in (\ref{eqn2034}). By (\ref{eqn2038}) we have
$\overline{T}=\overline{T}S_{\sigma}$ for any $\sigma\in{\cal B}$
and hence $\overline{T}=n^{-1}1_p1_t'$. Hence we have $\hat{\overline
{T}}=\overline{T}B_t=0$. By the same argument we have $\hat{\overline
{F}}=0$. Then the equality in (\ref{eqn2032}) holds, and so does the
equality in (\ref{eqn2036}). Hence we proved $\operatorname{
tr}(C_d)=q_d^*$.\vadjust{\goodbreak}

Given any design $d$ with corresponding $P=(p_s,s\in{\cal S})$, we
could define a new design $\bar{d}\leftrightarrow P_{\bar{d}}=(\bar
{p}_s,s\in{\cal S})$ by
%
\begin{equation}
\label{eqn2122} P_{\bar{d}}=\frac{\sum_{\sigma\in{\cal P}} P_{\sigma d}}{t!}.
\end{equation}
Then we have $\sum_{s\in{\cal S}}p_sq_s(x)=\sum_{s\in{\cal S}}\bar
{p}_sq_s(x)$ in view of (\ref{eqn2042}) and $q_d^*=q^*_{\bar{d}}$.
\end{pf*}

\begin{pf*}{Proof of Theorem~\ref{thm2042}} In the following, we would
apply Lemma 3.1 of \citet{Kus97N1} to prove (iii). The proof of (i)
and (ii) follows from similar arguments. Given any sequence $s$, we
have $q_s(x)=\sum^p_{k=m}\alpha_k q_s^k(x)$ where
$q_s^k(x)=q_{s11}^k+2q_{s12}^kx+q_{s22}^kx^2$ and
$q_{sij}^k=\operatorname{ tr}(G_i'B^k_pG_j)$ with $G_1=\hat{T}_u$ and $G_2=\hat
{F}_u$. By direct calculation we have
\begin{eqnarray*}
q_{s11}^k&=&k-\xi_{s_k}/k,
\\
q_{s12}^k&=&(k\rho_{s_k}+f_{s_k,t_k}-
\xi_{s_k})/k,
\\
q_{s22}^k&=&(kt-1) (k-1)/kt-(\xi_{s_k}-2f_{{s_k},t_k}+1)/k,
\end{eqnarray*}
where $\xi_{s_k}=\sum^t_{i=1}(f_{s_k,i})^2$ and $\rho_{s_k}=\sum^{k-1}_{j=1}1_{t_j=t_{j+1}}$.
For notational simplicity we define $\xi_k=\xi_{s_k}$, $\rho_k=\rho_{s_k}$ and $f_k=f_{s_k,t_k}$. Also let
$\Xi_A, A\subset\{x,k,p,t\}$, denote a quantity that depends on the
elements of $A$, and $a\propto_k b$ means that $a/b$ is a quantity
that only depend on $k$. Then
%
\begin{eqnarray}
q_s^k(x)&\propto_k&-\xi_k(x+1)^2+2f_k
\bigl(x+x^2\bigr)+2k\rho_kx+\Xi_{k,t,x}
\label{eqn207}
\\
&=&-(\xi_k-2f_k) (x+1)^2+2k(
\rho_k-f_k)x
\nonumber
\\
&&{}+2f_k\bigl[(k-1)x-1\bigr]+\Xi_{k,t,x}.\label{eqn2072}
\end{eqnarray}
From (\ref{eqn207}), for any $x>0$, the sequence which maximizes
$q_s^k(x)$ has to be of the form $
(1*1_{f_{s_k,1}}|2*1_{f_{s_k,2}}|\ldots\ldots|(t-1)*1_{f_{s_k,t-1}}|t*1_{f_{s_k,t}}
)$ with\vspace*{1.5pt} the restrictions of $f_{s_k,i+1}\geq f_{s_k,i},i=1,2,\ldots,t-1$
and $f_{s_k,t-1}-f_{s_k,1}\leq1$. For the special case of $k\leq t$,
the sequence reduces to the form of $(1,2,\ldots,k-h,t*1_h).$ By (\ref
{eqn2072}) the sequence of $\langle re\rangle$ maximizes $q_s(x)$ for
any $x\in((p-1)^{-1},(p-2)^{-1}]$ since this sequence maximizes
$q_s^k(x)$ for all $k=m,\ldots,p$. Since all the sequences in the class of
$\langle re\rangle$ have the same value of $d(q_s(x))/dx$, we need to
choose $x$ so that the derivative is zero, and hence (iii) is proven.
\end{pf*}

\begin{pf*}{Proof of Theorem~\ref{thm2142}} By Lemma~\ref{lemma215},
equations (\ref{eqn11})--(\ref{eqn44}) is equivalent to
%
\begin{eqnarray}
\check{C}_{d11}+x^*\check{C}_{d12}B_t&=&
\frac{ny^*}{t-1}B_t,\label{eqn1}
\\
\check{C}_{d21}+x^*\check{C}_{d22}B_t&=&0,
\label{eqn2}
\\
B\bigl(\hat{\overline{T}}+x^*\hat{\overline{F}}\bigr)&=&0.\label{eqn3}
\end{eqnarray}

First we show the \textit{necessity}. Let $f$ be a symmetric optimal
design and $g$ be a new design with $P_g=P_d/2+P_f/2$. Then by Lemmas
\ref{lemma2122} and~\ref{lemma212} we\break have
%
\begin{eqnarray}\label{eqn21211}
\check{C}_g&\geq& \check{C}_d/2+\check{C}_f/2
\nonumber
\\[-8pt]
\\[-8pt]
\nonumber
&=&\frac{ny^*}{t-1}B_t.
\end{eqnarray}
Let $\bar{g}$ be the symmetrized version of design $g$ as defined by
(\ref{eqn2122}). Following the same argument as in Lemma \ref
{lemma2122} we have
%
\begin{equation}
\label{eqn21212} \sum_{\sigma\in{\cal P}}S_{\sigma}'
\check{C}_gS_{\sigma}/|{\cal P}|\leq C_{\bar{g}}.
\end{equation}
Combining (\ref{eqn21211}) and (\ref{eqn21212}) we have
\[
C_{\bar{g}}=\frac{ny^*}{t-1}B_t,
\]
in view of Corollary~\ref{cor204}(ii). Then we have $\operatorname{ tr}(\check
{C}_d)=ny^*$ which together with (\ref{eqn21211}) yields
\[
\check{C}_g=\frac{ny^*}{t-1}B_t.
\]
Following similar arguments as in Theorem 5.3 of \citet{Kus97N2} we have
%
\begin{eqnarray}
\check{C}_{f22}\check{C}_{g22}^+\check{C}_{g21}&=&
\check {C}_{f21},\label{eqn21213}
\\
\check{C}_{d22}\check{C}_{g22}^+\check{C}_{g21}&=&
\check {C}_{d21},\label{eqn21214}
\end{eqnarray}
where $G^+$ denotes the Moore--Penrose inverse of $G$. Since $f$ is a
symmetric design, we have $\check{C}_{f21}=q_{f12}B_t/(t-1)$ and
$\check{C}_{f22}=q_{f22}B_t/(t-1)+(1_t'\check{C}_{f22}1_t)J_t/t^2$.
So we have $\check{C}_{f22}^+=(t-1)B_t/q_{f22}+J_t/(1_t'\check
{C}_{f22}1_t)$. By left multiplying both sides of (\ref{eqn21213}) we have
%
\begin{eqnarray}\label{eqn21215}
\check{C}_{g22}^+\check{C}_{g21}&=&\check{C}_{f22}^+
\check {C}_{f21}
\nonumber
\\[-8pt]
\\[-8pt]
\nonumber
&=&-x^*B_t.
\end{eqnarray}
By plugging (\ref{eqn21215}) into (\ref{eqn21214}) we have (\ref
{eqn2}). Then we have
\begin{eqnarray*}
\frac{ny^*}{t-1}B_t=\check{C}_{d}&=&
\check{C}_{d11}-\check {C}_{d12}\check{C}_{d22}^-
\check{C}_{d21}
\\
&=&\check{C}_{d11}+x^*\check{C}_{d12}\check{C}_{d22}^-
\check {C}_{d22}B_t
\\
&=&\check{C}_{d11}+x^*\check{C}_{d12}B_t.
\end{eqnarray*}
Hence (\ref{eqn1}) is derived. From (\ref{eqn2034}) and (5.3) of
\citet{Kus97N2} we have
%
\begin{eqnarray}\label{eqn2102}
ny^*&=&\operatorname{ tr}(\overline{C}_d)
\nonumber
\\
&\leq& \operatorname{ tr}(\tilde{C}_d)
\nonumber
\\[-8pt]
\\[-8pt]
\nonumber
&\leq&\operatorname{ tr}\bigl(\overline{C}_{d11}+2x\overline{C}_{d12}+x^2
\overline {C}_{d22}\bigr)
\nonumber
\\
&=&q_d(x)-n\operatorname{tr}\bigl[(\hat{\overline{T}}_d+x\hat{
\overline{F}}_d)'B(\hat {\overline{T}}_d+x
\hat{\overline{F}}_d)\bigr].\nonumber
\end{eqnarray}
Setting $x=x^*$ in (\ref{eqn2102}) gives $y^*\leq y^*-\operatorname{ tr}[(\hat
{\overline{T}}_d+x^*\hat{\overline{F}}_d)'B(\hat{\overline
{T}}_d+x^*\hat{\overline{F}}_d)]$ which yields (\ref{eqn3}) due to
Pukelsheim [(\citeyear{Puk93}), page 15].

Now we show the \textit{sufficiency}. By utilizing (\ref{eqn1}), (\ref
{eqn2}) and (\ref{eqn3}) we have
\begin{eqnarray*}
C_{d11}+x^*C_{d12}B_t&=&\frac{ny^*}{t-1}B_t,
\\
C_{d21}+x^*C_{d22}B_t&=&0,
\end{eqnarray*}
which in turn yields
\begin{eqnarray*}
C_d&=&C_{11}+x^*C_{d12}C_{d22}^-C_{d22}B_t
\\
&=&\frac{ny^*}{t-1}B_t.
\end{eqnarray*}
\upqed\end{pf*}

\section*{Acknowledgements}
We are grateful to the referees and the Associate Editor for their constructive comments on earlier
versions of this manuscript.

%


%
%

\printaddresses

\end{document}